\documentclass[a4paper,12pt]{amsart}
\usepackage{amssymb}
\usepackage{graphicx}
\usepackage{hyperref}
\usepackage{amsmath}
\usepackage{amsthm}
\usepackage{amssymb}
\newtheorem{definition}{Definition}[section] 
\hypersetup{colorlinks =true,citecolor=blue,linkcolor=blue,filecolor=blue,urlcolor=blue,citecolor=blue}
\usepackage{ifthen}
\nonstopmode \numberwithin{equation}{section}
\newtheorem{thm}{Theorem}
\newtheorem{cor}{Corollary}
\newtheorem{lem}{Lemma}

\newtheorem{conj}{Conjecture}

\theoremstyle{definition}

\newtheorem{prob}[equation]{Problem}

\newenvironment{rem}{%
	\bigskip
	\noindent \textsl{{\sl Remark. }}}{\bigskip}
\newenvironment{rems}{%
	\bigskip
	\noindent \textsl{{\sl Remarks. }}}{\bigskip}

\newcounter {own}
\def\theown {\thesection       .\arabic{own}}

\newenvironment{pf}[1][]{%
	\vskip 3mm
	\noindent
	\ifthenelse{\equal{#1}{}}%
	{{\slshape Proof. }}%
	{{\slshape #1.} }%
}%
{\qed\bigskip}

\newcounter{alphabet}

\newcommand{\A}{{\mathcal A}}

\newcommand{\ID}{{\mathbb D}}

\newcommand{\IC}{{\mathbb C}}

\def\be{\begin{equation}}
	\def\ee{\end{equation}}

\newcommand{\bee}{\begin{enumerate}}
	\newcommand{\eee}{\end{enumerate}}

\newcommand{\blem}{\begin{lem}}
	\newcommand{\elem}{\end{lem}}
\newcommand{\bthm}{\begin{thm}}
	\newcommand{\ethm}{\end{thm}}
\newcommand{\bcor}{\begin{cor}}
	\newcommand{\ecor}{\end{cor}}
\newcommand{\beg}{\begin{examp}}
	\newcommand{\eeg}{\end{examp}}
\newcommand{\begs}{\begin{examples}}
	\newcommand{\eegs}{\end{examples}}
\newcommand{\bdefe}{\begin{defin}}
	\newcommand{\edefe}{\end{defin}}
\newcommand{\bprob}{\begin{prob}}
	\newcommand{\eprob}{\end{prob}}
\newcommand{\bei}{\begin{itemize}}
	\newcommand{\eei}{\end{itemize}}

\newcommand{\bcon}{\begin{conj}}
	\newcommand{\econ}{\end{conj}}
\newcommand{\bcons}{\begin{conjs}}
	\newcommand{\econs}{\end{conjs}}
\newcommand{\bprop}{\begin{propo}}
	\newcommand{\eprop}{\end{propo}}
\newcommand{\br}{\begin{rem}}
	\newcommand{\er}{\end{rem}}
\newcommand{\brs}{\begin{rems}}
	\newcommand{\ers}{\end{rems}}
\newcommand{\bo}{\begin{obser}}
	\newcommand{\eo}{\end{obser}}
\newcommand{\bos}{\begin{obsers}}
	\newcommand{\eos}{\end{obsers}}
\newcommand{\bpf}{\begin{pf}}
	\newcommand{\epf}{\end{pf}}
\newcommand{\ba}{\begin{array}}
	\newcommand{\ea}{\end{array}}
\newcommand{\beq}{\begin{eqnarray}}
	\newcommand{\beqq}{\begin{eqnarray*}}
		\newcommand{\eeq}{\end{eqnarray}}
	\newcommand{\eeqq}{\end{eqnarray*}}

\newcounter{minutes}
\divide\time by 60
\newcounter{hours}
\multiply\time by 60 \addtocounter{minutes}{-\time}

\begin{document}
\bibliographystyle{amsplain}
\title{Landau-type Theorem and Bloch constant for elliptic harmonic mappings}
\begin{center}
{\tiny \texttt{FILE:~\jobname .tex,
printed: \number\year-\number\month-\number\day,
\thehours.\ifnum\theminutes<10{0}\fi\theminutes}
}
\end{center}
	
\begin{center}
		
\end{center}
\author{Bappaditya Bhowmik}
\address{Bappaditya Bhowmik, Department of Mathematics,
Indian Institute of Technology Kharagpur, Kharagpur-721302}
\email{bappaditya@maths.iitkgp.ac.in}
\author{Sambhunath Sen}
\address{Sambhunath Sen, Department of Mathematics,
Indian Institute of Technology Kharagpur, Kharagpur-721302}
\email{sensambhunath4@iitkgp.ac.in}
	

\begin{abstract}
In this article, we  obtain certain  estimates for the Taylor coefficients of $(K,K')$-elliptic harmonic mappings and using these estimates, we prove a Landau-type theorem for these mappings. We also derive Bloch constant for the class of $(K, K')$-elliptic harmonic mappings.
\end{abstract}
\thanks{The first author of this article would like to thank SERB, India for its financial support through Core Research Grant (Ref.No.-  CRG/2022/001835). The second author would like to thank the financial support from CSIR, HRDG, India (Ref. No.- 09/081(1389)/2019-EMR-I)}
	
\maketitle
\pagestyle{myheadings}
\markboth{B. Bhowmik and S. Sen}{Landau-type Theorem and Bloch constant for elliptic harmonic mappings}
	
\bigskip
 
\section{Introduction}
Let $\IC$ be the finite complex plane. By $\ID$, we denote the open unit disk of the complex plane. Let $\partial\ID$ be the unit circle and  $\ID_r:=\{z\in\IC :|z|<r\}$. In geometric function theory, a long-standing open problem is finding the exact value of the largest radius of  univalent disk contained in the image of a holomorphic function. A disk $\Omega$ will be called a univalent disk in $f(\ID)$ if there exists a domain $\Delta$  in $\ID$ such that $f$ maps univalently $\Delta$ onto $\Omega$. 
In 1924, André Bloch  proved that, if $f$ is a holomorphic function in $\ID$ with $f'(0)\neq 0$, then $f(\ID)$ contains a univalent disk of positive radius (c.f. \cite{Bloch}) ). This result is popularly known as the {\it Bloch Theorem}. The set of all analytic functions in a domain $D\subseteq\IC$ will be denoted by $\A(D)$.
Let $B_f$ be the largest radius of the univalent disk contained in $f(\ID)$ and
$$
B:=\sup\,\{B_f: f\in\A(\ID),\; f'(0)=1\}.
$$
We call $B$ as the Bloch constant for functions in the class $\mathcal{F}:=\{f\in \A(\ID): f'(0)=1\}$.
The latest known lower and upper bounds  for $B$ are  
$$
\frac{\sqrt{3}}{4}+2\times 10^{-4}< B \leq \frac{1}{\sqrt{1+\sqrt{3}}}\frac{\Gamma(1/3)\Gamma(11/12)}{\Gamma(1/4)}\approx 0.4719,
$$
which were obtained by Chen and Gauthier (c.f. \cite{CG}) and Ahlfors and Grunsky (c.f. \cite{AG}), respectively. Let $f\in\A(\ID)$ with $|f(z)|\leq M$ for all $z\in\ID$ and $f(0)=f'(0)-1=0$, then the classical \textit{Landau Theorem} asserts that $f$ is univalent in $\ID_{r_0}$ and $f(\ID_{r_0})$ contains the disk $\ID_{R_0}$, where, 
$$
r_0=\frac{1}{M+\sqrt{M^2-1}},\;\,\mbox{and}\;\, R_0=Mr_0^2.
$$
These two radii are best possible as can be seen by considering
$$
f(z)=Mz\left(\frac{1-Mz}{M-z}\right),\;z\in\ID.
$$

In this article, we obtain certain Taylor coefficient estimates for functions in the class of  $(K,K')$- elliptic harmonic mappings, and then using this estimate, we prove a Landau-type theorem for this class of mappings. We mention here that Theorems~\ref{thm1} and \ref{thm2} of this article improve  two recent results obtained by  Kumar et al. (c.f., \cite[Theorems~2.1, and 2.2]{AK}) for the class of $(K,K')$-elliptic harmonic mappings. More importantly, we also obtain Bloch-type theorems for $(K,K')$-elliptic harmonic mappings and in particular, these theorems generalize   certain results of M.-S. Liu (c.f. \cite[Theorems~2.4 and 2.5]{Liu}) for $K$-quasiregular harmonic mappings. 

 In order to define $(K,K')$-elliptic harmonic mappings, we need to prepare some backgrounds.
For a real $2\times 2$ matrix $A$, the matrix norm is defined by 
$$
\|A\|:=\sup\{|Az|:|z|=1\}.
$$
Next, the formal derivative of a complex-valued function $f=u+iv$ is  
$$
D_f=\begin{pmatrix}
	u_x & u_y\\
	v_x & v_y
\end{pmatrix}.
$$
For a continuously differentiable function $f$, we know that its partial derivatives can be computed as follows:
$$
f_z=\frac{1}{2}(f_x-if_y),\;\; f_{\overline{z}} =\frac{1}{2}(f_x+if_y),
$$
where, $ z=x+iy$. The Jacobian of $f$ can be computed as $J_f:=|f_z|^2-|f_{\overline{z}}|^2$.
\begin{definition}
	A sense-preserving and continuous mapping $f$ from $\ID$ onto $\IC$ is said to be $(K,K')$-elliptic mapping if
	\begin{itemize}
		\item[(i)] $f$ is absolutely continuous on lines in $\ID$, i.e., for every closed rectangle $R\subset\ID$, $f$ is  absolutely continuous on almost every line segment in $R$, which is parallel to the co-ordinate axes,
		\item[(ii)] $f_z, f_{\bar z}\in L^2(\ID)$ and $J_f\neq 0$ a.e. in $\ID$,
		\item[(iii)] there exist $K\geq 1$ and $K'\geq 0$ such that
		$$
		\|D_f\|^2\leq KJ_f+K'\;\;\mbox{a.e. in}\;\ID.   
		$$
	\end{itemize}
	If $K'=0$, then $f$ is called $K$-quasiregular mapping in $\ID$.
\end{definition}

We refer to the articles \cite{CP,FS,N} and references therein for more details about the $(K, K')$-elliptic mappings. 
We now turn our attention to harmonic mappings. 
A mapping $f$ is said to be a harmonic mapping on a domain $D\subseteq\IC$ if $f$ is twice continuously differentiable and satisfies $\bigtriangleup f=4f_{z\overline{z}}=0 $ on $D$. 
For a continuously differentiable function $f$, we define the maximum and minimum length distortions as
\beq\nonumber
\Lambda_f(z)&:=&\max_{0\leq \theta\leq 2\pi}\left|f_z(z)+e^{-2i\theta}f_{\overline{z}}(z)\right|=\left|f_z(z)\right|+\left|f_{\overline{z}}(z)\right|,\;\mbox{and}\\\nonumber
\lambda_f(z)&:=&\min_{0\leq \theta\leq 2\pi}\left|f_z(z)+e^{-2i\theta}f_{\overline{z}}(z)\right|=\left||f_z(z)|-|f_{\overline{z}}(z)|\right|.
\eeq
If a harmonic mapping is $(K, K')$-elliptic, then it is called $(K, K')$-elliptic harmonic mapping. 
A harmonic mapping $f$ is said to be $K$-quasiregular $(K\geq 1)$ on a domain $D$ if $\Lambda_f(z)\leq K\lambda_f(z)$  holds everywhere on $D$. 
Let $\mathcal{H}(\ID)$ is the set of all harmonic mappings $f=h+\overline{g}$ on $\ID$, where, $h,g\in\A(\ID)$  with the Taylor's series expansions
\be\label{eq1.1}
h(z)=\sum_{n=0}^\infty a_nz^n,\;\mbox{and}\; g(z)=\sum_{n=1}^\infty b_nz^n,\;z\in\ID.
\ee
A mapping $f$ defined on a domain $D\subseteq\IC$ is said to be an open mapping if $f$ maps any open subset of $D$ to an open set in $\IC$. It is well-known that, every non-constant holomorphic function is an open map, but for harmonic mappings, this
result is not always true. In 2000, Chen, Gauthier, and Hengartner (c.f., \cite{CGH})  obtained Landau-type theorems and Bloch constant for the open planar harmonic mappings and $K$-quasiregular harmonic mappings in $\ID$.
In 2009, Ming-Sheng Liu (c.f., \cite{Liu}) improved the results obtained by Chen et al. (c.f., \cite{CGH}).
We refer to \cite{BL, BS-1, BS-2, CPW, AK, LP, Liu2} and references therein on recent investigations on Bloch-type theorems and Landau-type theorems for different classes of mappings such as meromorphic, bi-analytic, polyharmonic, quasiregular, etc.. 

%
%
%
%
%
We organize this article as follows. In Section~\ref{sec2}, we present  main results of this article. In Section~\ref{sec3}, we prove all the main results.
\section{Main Results}\label{sec2}
\bthm\label{thm1}
Let $f=h+\overline{g}$ be a $(K,K')$-elliptic harmonic mapping in $\ID$ with $f(0)=\lambda_f(0)-1=0$. Let $h$ and $g$ have the Taylor's series expansion of the form $(\ref{eq1.1})$. If $\lambda_f(z)\leq\lambda$ for all $z\in\ID$, then 
$$
|a_n|+|b_n|\leq \frac{1}{n}\left(K\lambda+\frac{2(K'-1)}{K\lambda+\sqrt{K^2\lambda^2+4K'}}\right),\;n=2,3,\dots.
$$
Equality holds in the above inequality  for $K=1$, $K'=0$ and $n=2,3,\dots$ for the function
$$
F_n(z)=\lambda^2z-\int_0^z\frac{(\lambda^3-\lambda)dz}{\lambda+z^{n-1}},\; z\in\ID.
$$
\ethm
\br
For all $K\geq 1$, $K'\geq 0$ and $\lambda>1$, it is easy to show that
$$
\frac{2(K'-1)}{K\lambda+\sqrt{K^2\lambda^2+4K'}}<\sqrt{K'}.
$$
Thus, by virtue of the Theorem~\ref{thm1}, we get improved  coefficient estimates for $(K,K')$-elliptic harmonic mappings that was obtained in \cite[Theorem~2.1]{AK}.
\er

Using Theorem~\ref{thm1}, we get the following Landau-type theorem for the class of $(K,K')$-elliptic harmonic mappings.
\bthm\label{thm2}
Let $f=h+\overline{g}$ be a $(K,K')$-elliptic harmonic mapping in $\ID$ with $f(0)=0$, $\lambda_f(0)=1$ and $\lambda_f(z)\leq\lambda$ for all $z\in\ID$. Then $f$ is univalent in the disk $\ID_{r_1}$, where 
$$
r_1=\frac{1}{1+ K\lambda+\frac{2(K'-1)}{K\lambda+\sqrt{K^2\lambda^2+4K'}}},
$$
and $f(\ID_{r_1})$ contains a disk $\ID_{\sigma_1}$, where,
$$
\sigma_1= 1+ \left(K\lambda+\frac{2(K'-1)}{K\lambda+\sqrt{K^2\lambda^2+4K'}}\right)\ln\left( K\lambda r_1+\frac{2(K'-1)r_1}{K\lambda+\sqrt{K^2\lambda^2+4K'}}\right).
$$
\ethm

\brs
(i) Theorem~\ref{thm2}  improves the Landau-type theorem for the class of $(K,K')$-elliptic harmonic mappings in $\ID$ obtained in \cite[Theorem~2.2]{AK}, which we conclude by the following observation. It is easy to show that, for all $K\geq1$, $K'\geq 0$ and $\lambda>1$,
$$
0<\frac{K\lambda+\sqrt{K^2\lambda^2+4K'}}{2}\leq K\lambda+\sqrt{K'}.
$$
Thus, from the above inequality we get
\beq\nonumber\label{eq1.3}
K\lambda+\sqrt{K'}&>&K\lambda+\sqrt{K'}-\frac{2}{K\lambda+\sqrt{K^2\lambda^2+4K'}}\\\nonumber
&\geq& \frac{K\lambda+\sqrt{K^2\lambda^2+4K'}}{2}-\frac{2}{K\lambda+\sqrt{K^2\lambda^2+4K'}}\\
&=& K\lambda+\frac{2(K'-1)}{K\lambda+\sqrt{K^2\lambda^2+4K'}}.
\eeq
Therefore, the inequality (\ref{eq1.3}) implies that 
$$
r_1>\frac{1}{K\lambda+\sqrt{K'}}.
$$

Now to prove 
$$
\sigma_1>1+ \left(K\lambda+\sqrt{K'}\right)\ln\left(\frac{K\lambda+\sqrt{K'}}{1+K\lambda+\sqrt{K'}}\right)=\sigma\;\mbox{(say)},
$$
we consider the function
$$
\psi(t)=1+t\ln\left(\frac{t}{1+t}\right),\; t>0.
$$
It is easy to prove that $\psi'(t)<0$ for all $t>0$. Therefore, $\psi(t)$ is a decreasing function for $t>0$. Thus, from the inequality (\ref{eq1.3}), we get $\sigma_1>\sigma$.\\
(ii) If we take $K'=0$ in Theorem~\ref{thm2}, we get back the result proved by Liu (c.f., \cite[Theorem~2.2]{Liu}) for $K$-quasiregular harmonic mappings.
\ers

Now, we state the Bloch-type theorem for the class of $(K,K')$-elliptic harmonic mappings in $\ID$.

\bthm\label{thm3}
Let $f=h+\overline{g}$ be a $(K,K')$-elliptic harmonic mapping in $\ID$ with $f(0)=0$, $\lambda_f(0)=1$. Then $f(\ID)$ contains a univalent disk of radius at least
$$
\rho_1=\frac{1}{\sqrt{2}}+\left(\sqrt{2}K+\frac{4K'-1}{\sqrt{2}\left(K+\sqrt{K^2+4K'}\right)}\right)\ln\left(\frac{2K+\frac{4K'-1}{K+\sqrt{K^2+4K'}}}{1+2K+\frac{4K'-1}{K+\sqrt{K^2+4K'}}}\right).
$$
\ethm

\brs
(i) From this theorem, we conclude that, the lower bound of the Bloch constant for the class of $(K,K')$-elliptic harmonic mapping is $\rho_1$.\\
(ii) If we take $K'=0$ in the above theorem, we get back the lower bound of the Bloch constant obtained by Liu (c.f., \cite[Theorem~2.4]{Liu}) for $K$-quasiregular harmonic mappings.
\ers
\bcor\label{cor1}
Let $f=h+\overline{g}$ be a $(K,K')$-elliptic harmonic mapping in $\ID$ with $f(0)=0$, $J_f(0)=1$. Then $f(\ID)$ contains a univalent disk of radius at least 
$$
\rho_2=\frac{1}{\sqrt{2(K+K')}}\left(1+ t(K,K')\ln\left(\frac{t(K,K')}{1+t(K,K')}\right)\right),
$$
where,
$$
t(K,K')=\left(2K+\frac{4K'(K+K')-1}{K+\sqrt{K^2+4K'(K+K')}}\right).
$$
\ecor

\br
If we take $K'=0$ in the above theorem, we get back the lower bound of  Bloch constant obtained by Liu (c.f., \cite[Theorem~2.5]{Liu}) for the $K$-quasiregular harmonic mapping in $\ID$.
\er
\section{Proof of the main results}\label{sec3}
 In order to prove Theorem~\ref{thm1}, we will use the following Lemma obtained by Liu (c.f., \cite{Liu}).
\vspace{0.1cm}

\noindent\textbf{Lemma A.} (\cite[Lemma~2.1]{Liu}) \textit{Let $f=h+\overline{g}$ be a harmonic mapping in $\ID$ with $f(0)=\lambda_f(0)-1=0$. Let $h$ and $g$ have the Taylor's series expansion of the form $(\ref{eq1.1})$. If $\Lambda_f(z)\leq\Lambda$ for all $z\in\ID$, then
$$
|a_n|+|b_n|\leq\frac{\Lambda^2-1}{n\Lambda},\;\;n=2,3,\dots.
$$
The sharpness of the above estimates holds for the function
$$
f_n(z)=\Lambda^2z-\int_0^z\frac{(\Lambda^3-\Lambda)dz}{\Lambda+z^{n-1}},\; z\in\ID.
$$
}
\vspace{0.1cm}

\bpf[Proof of Theorem~\ref{thm1}]
Since $f$ is $(K,K')$-elliptic harmonic mapping, where $K\geq1$, $K'\geq0$ and $\lambda_f(0)=1$, therefore
$$
K\lambda\geq K\lambda(0)\geq 1.
$$
Since $\|D_f\|^2\leq KJ_f+K'$, we have 
$$
\Lambda_f(z)\leq \frac{1}{2}\left(K\lambda+\sqrt{K^2\lambda^2+4K'}\right)
$$
for all $z\in\ID$.
Therefore from the Lemma~A, we get
\be\label{eq3.1}
|a_n|+|b_n|\leq \frac{1}{n}\left(K\lambda+\frac{2(K'-1)}{K\lambda+\sqrt{K^2\lambda^2+4K'}}\right),\;n=2,3,\dots.
\ee
Now for $n=2,3,\dots$, the function
$$
F_n(z)=\lambda^2z-\int_0^z\frac{(\lambda^3-\lambda)dz}{\lambda+z^{n-1}}=z+\frac{\lambda^2-1}{n\lambda}z^n+\sum_{k=2}^\infty\frac{(-1)^{k+1}(\lambda^2-1)}{(kn-k+1)\lambda^k}z^{kn-k+1},\; z\in\ID,
$$
is $(1,0)$-elliptic harmonic mapping with $F_n(0)=0$, $\lambda_{F_n}(0)=1$ and $\lambda_{F_n}(z)\leq\lambda$ for all $z\in\ID$. Thus, the equality holds in (\ref{eq3.1}) for this function $F_n$, $n=2,3,\dots$.
This completes the proof.
\epf

\bpf[Proof of Theorem~\ref{thm2}]
Let us consider $A(z):=f(z)-f_z(0)z-f_{\overline{z}}(0)\overline{z}$, $z\in\ID$. Therefore, for all $z\in\ID$
$$
A(z)=\int_{[0,z]}(f_\zeta(\zeta)-f_\zeta(0))d\zeta+(f_{\overline{\zeta}}(\zeta)-f_{\overline{\zeta}}(0))d{\overline{\zeta}}.
$$
Now for $z_1,z_2\in\ID_r$, $0<r<1$,
\beq\nonumber
|A(z_1)-A(z_2)|&=&\left|\int_{[z_1,z_2]}(f_\zeta(\zeta)-f_\zeta(0))d\zeta+(f_{\overline{\zeta}}(\zeta)-f_{\overline{\zeta}}(0))d{\overline{\zeta}}\right|\\\nonumber
&=&\left|\int_{[z_1,z_2]}\left(\sum_{n=2}^\infty na_n\zeta^{n-1}\right)d\zeta+\left(\sum_{n=2}^\infty n\overline{b_n}\,\overline{\zeta}^{n-1}\right)d\overline{\zeta}\right|\\\nonumber
&\leq&\int_{[z_1,z_2]}\left(\sum_{n=2}^\infty n|a_n|r^{n-1}\right)|d\zeta|+\left(\sum_{n=2}^\infty n|b_n|r^{n-1}\right)|d\overline{\zeta}|\\\nonumber
&=&|z_1-z_2|\int_0^1\sum_{n=2}^\infty n(|a_n|+|b_n|)r^{n-1}dt.
\eeq
Applying Theorem~\ref{thm1} in the last inequality, we get
$$
|A(z_1)-A(z_2)|\leq |z_1-z_2|\left(K\lambda+\frac{2(K'-1)}{K\lambda+\sqrt{K^2\lambda^2+4K'}}\right)\frac{r}{1-r},
$$
for all $z_1,z_2\in\ID_r$, $0<r<1$.
Thus, for $z_1,z_2\in\ID_r$, using the last inequality, we obtain
\beq\nonumber
|f(z_1)-f(z_2)|&=&|A(z_1)-A(z_2)+(z_1-z_2)f_z(0)+\overline{(z_1-z_2)}f_{\overline{z}}(0))|\\\nonumber
&\geq& |z_1-z_2|\lambda_f(0)-|A(z_1)-A(z_2)|\\\nonumber
&\geq&|z_1-z_2|\left(1-\left(K\lambda+\frac{2(K'-1)}{K\lambda+\sqrt{K^2\lambda^2+4K'}}\right)\frac{r}{1-r}\right).
\eeq
Therefore, for $z_1\neq z_2$, $|f(z_1)-f(z_2)|>0$ in the disk $\ID_r$, whenever
$$
 \left(K\lambda+\frac{2(K'-1)}{K\lambda+\sqrt{K^2\lambda^2+4K'}}\right)\frac{r}{1-r} <1,
$$
then we have
$$
r<\frac{1}{1+ K\lambda+\frac{2(K'-1)}{K\lambda+\sqrt{K^2\lambda^2+4K'}}}=r_1.
$$
Hence, $f$ is univalent in $\ID_{r_1}$. This completes proof of the first part of the theorem.
Now to prove the second part of the theorem, we consider $z=r_1e^{i\theta}$. Then,
\beq\nonumber
|f(z)|&=&|A(z)+f_z(0)z+f_{\overline{z}}(0)\overline{z}|\\\nonumber
&\geq& \lambda_f(0)r_1-|A(z)|\\\nonumber
&=& r_1-\left|\int_{[0,z]}(f_\zeta(\zeta)-f_\zeta(0))d\zeta+(f_{\overline{\zeta}}(\zeta)-f_{\overline{\zeta}}(0))d{\overline{\zeta}}\right| \\\nonumber
&=&r_1-\int_0^1\sum_{n=2}^\infty n(|a_n|+|b_n|)|zt|^{n-1}|z|dt.
\eeq
Applying Theorem~\ref{thm1} in the last inequality, we get for $z=r_1e^{i\theta}$,
\beq\nonumber
|f(z)|&\geq& r_1-\left(K\lambda+\frac{2(K'-1)}{K\lambda+\sqrt{K^2\lambda^2+4K'}}\right)\int_0^1\frac{|z|^2tdt}{1-|z|t}\\\nonumber
&=& 1+ \left(K\lambda+\frac{2(K'-1)}{K\lambda+\sqrt{K^2\lambda^2+4K'}}\right)\ln\left( K\lambda r_1+\frac{2(K'-1)r_1}{K\lambda+\sqrt{K^2\lambda^2+4K'}}\right)\\\nonumber
&=&\sigma_1.
\eeq
Therefore,  
$f(\ID_{r_1})$ contains the disk $\ID_{\sigma_1}$. This completes the proof.
\epf

\bpf[Proof of Theorem~\ref{thm3}]
Without loss of generality, we assume that $f$ is harmonic on $\partial\ID$ and 
$$
M=\sup_{z\in\ID}(1-|z|^2)\lambda_f(z).
$$
Since, $\lambda_f(z)$ is continuous on $\overline{\ID}$, then there exists a point $z_0\in\ID$ such that 
$$
M=(1-|z_0|^2)\lambda_f(z_0).
$$
Let, $\phi(z)=(z+z_0)/(1+\overline{z_0}z)$, $z\in\ID$, then $\phi(0)=z_0$. Let
$$
F(\zeta):=\frac{1}{M}(f(\phi(\zeta))-f(z_0)),~ \zeta\in\ID.
$$
Thus, $F$ is harmonic in $\ID$, $F(0)=0$ and $F(\ID)=f(\ID)/M$. Also, $\lambda_F(0)=1$ and 
$$
(1-|\zeta|^2)\lambda_F(\zeta)=\frac{1}{M}(1-|\phi(\zeta)|^2)\lambda_f(\phi(\zeta))\leq1.
$$
Let $G(w):=\sqrt{2}F(w/\sqrt{2})$, $w\in\ID$. Then $G(0)=0$ and for $w\in\ID$,
\beq\nonumber
\Lambda_G(w)&=&\Lambda_F\left(\frac{w}{\sqrt{2}}\right)=\frac{1}{M}\left|\phi'\left(\frac{w}{\sqrt{2}}\right)\right|\Lambda_f\left(\phi\left(\frac{w}{\sqrt{2}}\right)\right),\;\mbox{and}\\\nonumber
\lambda_G(w)&=&\lambda_F\left(\frac{w}{\sqrt{2}}\right)=\frac{1}{M}\left|\phi'\left(\frac{w}{\sqrt{2}}\right)\right|\lambda_f\left(\phi\left(\frac{w}{\sqrt{2}}\right)\right).
\eeq
Now for $w\in\ID$,
\beq\nonumber
\Lambda_G^2(w)&=&\frac{1}{M^2}\left|\phi'\left(\frac{w}{\sqrt{2}}\right)\right|^2\Lambda_f^2\left(\phi\left(\frac{w}{\sqrt{2}}\right)\right)\\\nonumber
&\leq& \frac{1}{M^2}\left|\phi'\left(\frac{w}{\sqrt{2}}\right)\right|^2\left(KJ_f\left(\phi\left(\frac{w}{\sqrt{2}}\right)\right)+K'\right)\\\nonumber
&=& KJ_G(w)+\frac{K'}{M^2}\left|\phi'\left(\frac{w}{\sqrt{2}}\right)\right|^2\\\nonumber
&\leq& KJ_G(w)+4K'.
\eeq
Thus, $G$ is $(K,4K')$-elliptic harmonic mapping in $\ID$ with $\lambda_G(0)=\lambda_F(0)=1$, $G(0)=0$ and for all $z\in\ID$,
$$
\lambda_G(z)\leq\frac{2}{2-|z|^2}\leq 2.
$$
Therefore, using Theorem~\ref{thm2} we conclude that $G(\ID)$ contains a univalent disk of radius 
$$
R=1+\left(2K+\frac{4K'-1}{K+\sqrt{K^2+4K'}}\right)\ln\left(\frac{2K+\frac{4K'-1}{K+\sqrt{K^2+4K'}}}{1+2K+\frac{4K'-1}{K+\sqrt{K^2+4K'}}}\right).
$$
As $f(\ID)=MF(\ID)$ and $G(\ID)=\sqrt{2}F(\ID_{1/\sqrt{2}})$, thus, $f(\ID)$ contains a univalent disk of radius at least 
$$
\frac{MR}{\sqrt{2}}\geq \frac{\lambda_f(0)R}{\sqrt{2}}= \rho_1.
$$
This completes proof of the theorem.
\epf

\bpf[Proof of Corollary~\ref{cor1}]
Since, $f$ is a $(K,K')$-elliptic harmonic mapping with $f(0)=J_f(0)-1=0$, then
$$
\Lambda_f(0)\leq\sqrt{K+K'}.
$$
Thus, as $J_f(0)=1$, we get
$$
\lambda_f(0)=\frac{1}{\Lambda_f(0)}\geq \frac{1}{\sqrt{K+K'}}.
$$
Let
$$
F(z):=\frac{f(z)}{\lambda_f(0)},\;z\in\ID.
$$
Then $F(0)=\lambda_F(0)-1=0$ and $F$ is $(K,K'(K+K'))$-elliptic harmonic mapping in $\ID$.
Now applying Theorem~\ref{thm3} for the function $F$, we get the desired result. This completes the proof.
\epf

\noindent{\bf Declarations of interest}: None\\
\noindent{\bf Data availability statement}: Data sharing not applicable to this article as no datasets were generated or analysed during the current study.

\end{document}